\documentclass[a4paper,reqno]{amsart}
\usepackage{amssymb, amsmath, amscd}
\usepackage[final]{graphicx}
\usepackage{color}

\newtheorem{theorem}{Theorem}[section]

\newtheorem{example}[theorem]{Example}
\makeatletter
 \@addtoreset{equation}{section}
\makeatother
\begin{document}
\title{Anti-invariant Riemannian Submersions}
\author{P. Gilkey}
\address{PG: Mathematics Department, \; University of Oregon, \;\;
 Eugene \; OR 97403, \; USA}
\email{gilkey@uoregon.edu}
\author{M. Itoh}
\address{MI: University of Tsukuba,
1-1-1 Tennodai, Tsukuba, Ibaraki 305-8577 Japan}
\email{itohm@sakura.cc.tsukuba.ac.jp}
\author{J. H. Park}
\address{JHP:Department of
Mathematics, Sungkyunkwan University, Suwon, 440-746, Korea}
\email{parkj@skku.edu}
\begin{abstract}
We give a general Lie-theoretic
construction for anti-invariant almost Hermitian Riemannian submersions,
anti-invariant quaternion Riemannian submersions, anti-invariant para-Hermitian
Riemannian submersions, anti-invariant para-quaternion Riemannian submersions,
and anti-invariant octonian Riemannian submersions. This yields many compact Einstein examples.
\end{abstract}
\maketitle
\section{Introduction}
We begin by establishing some notational conventions.
\subsection{Riemannian submersions}
Let $M$ and $N$ be smooth manifolds of dimension $m$ and $n$, respectively,
and let $\pi:M\rightarrow N$ be a smooth map.
We say that $\pi$ is a {\it submersion} if $\pi_*$ is a surjective map from the tangent space
$T_PM$ to the tangent space $T_{\pi P}N$
for every point $P$ of $M$. Let $g_M$ and $g_N$ be Riemannian metrics on $M$ and $N$.
If $\pi:M\rightarrow N$ is a submersion, then the {\it vertical distribution} is the kernel of $\pi_*$
and the {\it horizontal distribution}
$\mathcal{H}$ is $\mathcal{V}^\perp$. We may then decompose $TM=\mathcal{V}\oplus\mathcal{H}$.
We say that $\pi$ is a {\it Riemannian submersion} if $\pi_*$ is an isometry from $\mathcal{H}_P$ to
$T_{\pi P}N$ for every point $P$ of $M$. We refer to O'Neill~\cite{O66} for further details
concerning the geometry of Riemannian submersions. If $g_M$ and $g_N$ are pseudo-Riemannian
metrics, we impose in addition the condition that the restriction of $g_M$ to $\mathcal{V}$ is non-degenerate
to ensure that $\mathcal{V}\cap\mathcal{H}=\{0\}$.
This gives rise to the notation of a
{\it pseudo-Riemannian} submersion.

\subsection{Hermitian geometry}
An endomorphism $J$ of $TM$
is said to define an {\it almost complex structure} on $M$
if $J^2=-\operatorname{id}$, i.e. $J$ gives a complex structure to $T_PM$ for
every point $P$ of $M$. We complexity the tangent bundle and let
$$
T^{1,0}:=\{X\in TM\otimes_{\mathbb{R}}\mathbb{C}:JX=\sqrt{-1}X\}\,.
$$
One says $J$ is {\it integrable} if $T^{1,0}$ is integrable, i.e.
$X,Y$ belong to $ C^\infty(T^{1,0})$ implies the complex Lie bracket
$[X,Y]$ also belongs to $C^\infty(T^{1,0})$. The Newlander-Nirenberg
Theorem \cite{NN57} is the analogue in the complex setting of the
Frobenius theorem in the real setting; $J$ is integrable if and only
if it arises from an underlying holomorphic structure on $N$. The
Riemannian metric $g_M$ is said to be {\it almost Hermitian} if
$J^*g_M=g_M$, i.e. if $g_M(JX,JY)=g_M(X,Y)$ for all tangent vectors
$X,Y\in T_PM$ and all points $P$ of $M$; the triple $(M,g_M,J)$ is
then said to be an {\it almost Hermitian manifold}; in the
pseudo-Riemannian setting one obtains the notion of {\it almost
pseudo-Hermitian} manifold similarly. The notation {\it Hermitian}
or {\it pseudo-Hermitian} is used if the structure $J$ is
integrable.

Let $(M,g_M,J)$ be an almost Hermitian manifold and let $\pi$ be a
Riemannian submersion from $(M,g_M)$ to $(N,g_N)$. Following the
seminal work of {\c{S}ahin}~\cite{S10,S13}, one says that $\pi$ is
an {\it anti-invariant almost Hermitian Riemannian submersion} if
$$J\{\mathcal{V}\}\subset\mathcal{H}\,.$$ If
$J\{\mathcal{V}\}=\mathcal{H}$, then $\pi$ is said to be {\it
Lagrangian}. There have been a number of subsequent papers in this
subject extending the work of {\c{S}ahin}~\cite{S10,S13}; we shall
cite just a few representative examples. Lee et al. \cite{LPSS15}
examined the geometry of anti-invariant Hermitian submersions from a
K\"ahler manifold onto a Riemannian manifold in relation to the
Einstein condition and examined when the submersions were Clairant
submersions. Ali and Fatima \cite{AF13} examined the nearly K\"ahler
setting. We also refer to related work of Ali and Fatima
\cite{AF13a}, of Beri et al.~\cite{BKEM15}, and of Murthan and
K\"upeli-Erken~\cite{MK13}.

\subsection{Quaternion geometry}
We shall restrict to flat quaternion structures as this is sufficient for our purposes.
 The {\it quaternion algebra} $\mathbb{Q}:=\mathbb{R}^4=\operatorname{Span}\{e_0,e_1,e_2,e_3\}$
 is defined by the relations:
\begin{equation}\label{Eq1.a}
\begin{array}{|r|r|r|r|r|}
\noalign{\hrule}
&e_0&e_1&e_2&e_3\\
\noalign{\hrule}e_0&e_0&e_1&e_2&e_3\\
\noalign{\hrule}e_1&e_1&-e_0&e_3&-e_2\\
\noalign{\hrule}e_2&e_2&-e_3&-e_0&e_1\\
\noalign{\hrule}e_3&e_3&e_2&-e_1&-e_0\\
\noalign{\hrule}
\end{array}\ .
\end{equation}
One says that $x\in\mathbb{Q}$ is an {\it imaginary quaternion} if $x\in\operatorname{Span}\{e_1,e_2,e_3\}$.
A {\it flat quaternion structure} on a manifold $M$ is a unital action of $\mathbb{Q}$ on $TM$.
If $x$ is a unit length purely imaginary quaternion, then $\xi\rightarrow x\cdot\xi$ defines an
almost complex structure on $M$.
If $g$ is a Riemannian metric on $M$, we shall assume in addition that $\|x\cdot\xi\|=\|x\|\cdot\|\xi\|$
for any quaternion $x$ and any tangent vector $\xi$.
Let $\pi:(M,g_M)\rightarrow(N,g_N)$ be a Riemannian submersion. Then one says
$\pi$ is an {\it anti-invariant quarternion Riemannian submersion} if
$x\cdot\mathcal{V}\subset\mathcal{H}$ for any purely imaginary quaternion $x$.
We have assumed that the roles of $\{e_1,e_2,e_3\}$ are globally defined (i.e. the structure
is flat); we refer to
Alekseevsky and Marchiafava~\cite{AM01} for a discussion of the more general setting.
Anti-invariant quaternion Riemannian submersions have been studied by K. Park \cite{P15}.

\subsection{Para-Hermitian geometry}
Instead of considering almost complex structures, one can consider
para-complex structures. Let $\tilde{\mathbb{C}}:=\mathbb{R}^2$ with
the para-complex structure $Je_1=e_2$ and $Je_2=e_1$. Let $(M,g_M)$
be a pseudo-Riemannian manifold of neutral signature $(\ell,\ell)$.
If $J$ is an endomorphism of $M$ with $J^2=\operatorname{Id}$ such
that $g_M(JX,JY)=-g_M(X,Y)$ for all $X,Y\in T_PM$ and all points $P$
of $M$, then the triple $(M,g_M,J)$ is said to be a {\it
para-Hermitian manifold}. Let $\pi$ be a pseudo-Riemannian
submersion from $(M,g_M)$ to $(N,g_N)$ with
$J\{\mathcal{V}\}\subset\mathcal{H}$. One then says $\pi$ is an {\it
anti-invariant para-Hermitian Riemannian submersion}; $\pi$ is {\it
Lagrangian para-Hermitian} if $J\{\mathcal{V}\}=\mathcal{H}$.
Atceken \cite{A13} and G\"und\"uzalp ~\cite{G13} examined this
setting. {G\"und\"uzalp~\cite{G14} also examined the anti-invariant
almost product setting; in the interests of brevity we shall not
treat this setting in this paper although our methods are clearly
applicable.}

\subsection{Para-quaternion geometry} In place of the quaternion
commutation relations given in Equation~(\ref{Eq1.a}), one imposes the para-quaternion relations
to define the para-quaternions
$\tilde{\mathbb{Q}}$ by setting:
$$
\begin{array}{|r|r|r|r|r|}
\noalign{\hrule}
&e_0&e_1&e_2&e_3\\
\noalign{\hrule}e_0&e_0&e_1&e_2&e_3\\
\noalign{\hrule}e_1&e_1&-e_0&e_3&-e_2\\
\noalign{\hrule}e_2&e_2&-e_3&{+}e_0&{ -}e_1\\
\noalign{\hrule}e_3&e_3&e_2&{+}e_1&{+}e_0\\
\noalign{\hrule}
\end{array}\ .
$$
If $J_1$ is Hermitian
and if $J_2$ and $J_3$ are para-Hermitian, then one obtains the
notion of a {\it para-quaternion manifold}. We refer to Ivanov and Zamkovoy \cite{IZ11}
for further details. If $\pi:(M,g_M)\rightarrow(N,g_N)$ is a pseudo-Riemannian submersion
and if $J_i\{\mathcal{V}\}\subset\mathcal{H}$ for $1\le i\le 3$, then $\pi$ is said to
be an {\it anti-invariant para-quaternion Riemannian submersion}.
To the best of our knowledge, there are no
papers on such geometries.

\subsection{Octonian geometry}
The octonians $\mathbb{O}$ arise from a non-associative and non-commutative bilinear multiplication
on $\mathbb{R}^8$.
If $\{e_0,\dots,e_7\}$ is the standard basis for $\mathbb{R}^8$, the multiplication is given
by the following table (see Wikipedia~\cite{W15}):
$$\begin{array}{|r|r|r|r|r|r|r|r|r|}
\noalign{\hrule}&e_0&e_1&e_2&e_3&e_4&e_5&e_6&e_7\\
\noalign{\hrule}e_0&e_0&e_1&e_2&e_3&e_4&e_5&e_6&e_7\\
\noalign{\hrule}e_1&e_1&-e_0&e_3&-e_2&e_5&-e_4&-e_7&e_6\\
\noalign{\hrule}e_2&e_2&-e_3&-e_0&e_1&e_6&e_7&-e_4&-e_5\\
\noalign{\hrule}e_3&e_3&e_2&-e_1&-e_0&e_7&-e_6&e_5&-e_4\\
\noalign{\hrule}e_4&e_4&-e_5&-e_6&-e_7&-e_0&e_1&e_2&e_3\\
\noalign{\hrule}e_5&e_5&e_4&-e_7&e_6&-e_1&-e_0&-e_3&e_2\\
\noalign{\hrule}e_6&e_6&e_7&e_4&-e_5&-e_2&e_3&-e_0&-e_1\\
\noalign{\hrule}e_7&e_7&-e_6&e_5&e_4&-e_3&-e_2&e_1&-e_0\\
\noalign{\hrule}
\end{array}\ .$$
The octonians satisfy the identity
$$
\|x\cdot y\|=\|x\|\cdot\|y\|\text{ for all }x,y\in\mathbb{R}^8\,.
$$
If $x\in\operatorname{Span}\{e_1,\dots,e_7\}$, then $x$ is said to be a {\it purely imaginary octonian}.
Such an octonian satisfies $x\cdot y\perp y$ for any $y\in\mathbb{R}^8$.
Let $(M,g)$ be a Riemannian manifold.
A {\it flat octonian} structure on a Riemannian manifold $(M,g)$ is
a unital octonian action on $TM$ such that $\|x\cdot\xi\|=\|x\|\cdot\|\xi\|$
for any octonian $x$ and any tangent vector $\xi$. If $\pi$
is a Riemannian submersion from $(M,g)$ to $(N,h)$,
 then we say that $\pi$ is {\it anti-invariant octonian} if $x\cdot\mathcal{V}\perp\mathcal{V}$
for any purely imaginary octonian $x\in\mathbb{R}^7$. To the best of our knowledge,
there are no papers dealing with anti-invariant octonian Riemannian submersions.

\subsection{Outline of the paper}
In Section~\ref{S2}, we will use Lie theoretic methods to construct
examples of anti-invariant almost Hermitian Riemannian submersions,
of anti-invariant quaternion Riemannian submersions, of
anti-invariant para-Hermitian Riemannian submersions, and of
anti-invariant para-quaternion Riemannian submersions. In
Section~\ref{S3}, we will discuss some examples which arise from
this construction. We conclude in Section~\ref{S4} by presenting a
different family of examples (including an anti-invariant octonian
Riemannian submersion) relating to the Hopf fibration where the
total space is not a Lie group. It is our hope that having a rich
family of examples will inform further investigations in this field.

\section{A Lie-theoretic construction}\label{S2}

Let $H$ be a closed and connected subgroup of an even dimensional
Lie group $G$. Let $\mathfrak{h}$ and $\mathfrak{g}$ be the
associated Lie algebras, respectively. Let
$\langle\cdot,\cdot\rangle=\langle\cdot,\cdot\rangle_{\mathfrak{g}}$
be a non-degenerate symmetric bilinear form {on} $\mathfrak{g}$
which is invariant under the adjoint action of $H$ and whose
restriction to $\mathfrak{h}$ is non-degenerate as well. We use
$\langle\cdot,\cdot\rangle_{\mathfrak{g}}$ to decompose
$\mathfrak{g}=\mathfrak{h}\oplus\mathfrak{h}^\perp$ as an orthogonal
direct sum. The inner product
$\langle\cdot,\cdot\rangle_{\mathfrak{g}}$ defines a left-invariant
pseudo-Riemannian metric on $G$ and, since the inner product is
invariant under the adjoint action of $H$, the restriction of
$\langle\cdot,\cdot\rangle_{\mathfrak{g}}$ to $\mathfrak{h}^\perp$
defines a $G$-invariant pseudo-Riemannian metric on the coset
manifold $G/H$ so that the natural projection $\pi:G\rightarrow G/H$
is a pseudo-Riemannian submersion.

\subsection{Complex geometry}
Assume $\langle\cdot,\cdot\rangle$ is
positive definite. Let
$J$ be a Hermitian complex structure on $\mathfrak{g}$; $J$ induces a left-invariant
Hermitian almost complex structure on $G$. Assume that
$J\{\mathfrak{h}\}\subset\mathfrak{h}^\perp$. Then $\pi:G\rightarrow G/H$ is
an anti-invariant almost Hermitian Riemannian submersion;
$\pi$ is Lagrangian if and only if $2\dim\{\mathfrak{h}\}=\dim\{\mathfrak{g}\}$.
More generally, if $\langle\cdot,\cdot\rangle$ is only assumed to be a non-degenerate
inner product and if the restriction to $\mathfrak{h}$ is assumed to be non-degenerate, then
we obtain an anti-invariant almost pseudo-Hermitian Riemannian submersion.

\subsection{Quaternion geometry}
Assume $\langle\cdot,\cdot\rangle$
is positive definite. Assume given a Hermitian quaternion structure
on $\mathfrak{g}$ such that $x\cdot\mathfrak{h}\subset\mathfrak{h}^\perp$ for any
purely imaginary quaternion $x$. Then
$\pi:G\rightarrow G/H$ is an anti-invariant quaternion Riemannian submersion.

\subsection{Para complex geometry}
Assume $\langle\cdot,\cdot\rangle$ has neutral signature and that
$\langle\cdot,\cdot\rangle$ is non-degenerate on $\mathfrak{h}$. Let $J$ be a Hermitian
para-complex structure on $\mathfrak{g}$ with $J\{\mathfrak{h}\}\subset\mathfrak{h}^\perp$.
Then $\pi:G\rightarrow G/H$
is an anti-invariant para-Hermitian Riemannian submersion; $\pi$ is Lagrangian if and only if
$2\dim\{\mathfrak{h}\}=\dim\{\mathfrak{g}\}$.

\subsection{Para-quaternion geometry} Assume $\langle\cdot,\cdot\rangle$ has neutral
signature and that $\langle\cdot,\cdot\rangle$ is non-degenerate on $\mathfrak{h}$.
Assume given a Hermitian para-quaternion structure on $\mathfrak{g}$ such
that $x\cdot\mathfrak{h}\subset\mathfrak{h}^\perp$ for any purely imaginary
para-quaternion $x$. Then $\pi:G\rightarrow G/H$
is an anti-invariant para-quaternion Riemannian submersion.

\section{Examples}\label{S3}
In this section, we present examples of anti-invariant Riemannian
submersions where the total space is a Lie group $G$ and the base
space is a homogeneous space upon which $G$ acts transitively by
isometries; $H\subset G$ is the isotropy subgroup of the action.
Example~\ref{Ex3.1} and Example~\ref{Ex3.2} are flat geometries.
Example~\ref{Ex3.3} arises from the Hopf fibration $S^1\rightarrow
S^3\rightarrow S^2$. In Example~\ref{Ex3.4}, the total space will be
$(S^3)^\nu$. We will take {product metrics} and if the metric on
$S^3$ is the usual round metric, these examples will be Einstein. In
Example~\ref{E3.6}, we take
$G=\mathbb{R}\times\operatorname{SL}(2,\mathbb{R})$ to construct
negative curvature examples.
\subsection{Abelian examples}
\begin{example}\label{Ex3.1}\rm Let $G=\mathbb{R}^m$ and let $H=\mathbb{R}^n\oplus0\subset G$
for $n<m$. We identify $G/H$ with $0\oplus\mathbb{R}^{m-n}$ and $\pi$ with
projection on the last $m-n$ coordinates.
\begin{enumerate}
\item Take the standard Euclidean inner product
on $G$ to obtain a bi-invariant Riemannian metric so that $\pi$ is a
Riemannian submersion.
\begin{enumerate}
\item Suppose $m=2\ell$
and $n=\ell$. Identify $G=\mathbb{C}^\ell$ so that $H$ corresponds to the
purely real vectors in $\mathbb{C}^\ell$. We identify $\mathfrak{g}$ with $G$
and $\mathfrak{h}$ with $H$. Then $\sqrt{-1}\mathfrak{h}\perp\mathfrak{h}$ and we
obtain a Lagrangian Hermitian Riemannian submersion; the almost complex structure
corresponds to scalar multiplication by $\sqrt{-1}$ and
is integrable.
\item Assume $m=4\ell$ and $n=\ell$. Identify $G=\mathbb{Q}^\ell$ so that
$H$ corresponds to the purely real vectors in $\mathbb{Q}^\ell$.
Then $x\cdot\mathfrak{h}\perp\mathfrak{h}$ if $x$ is a purely imaginary quaternion
and we obtain a Riemannian submersion which is anti-invariant quaternion.
\item Assume $m=8\ell$ and $n=\ell$. Identify $G=\mathbb{O}^\ell$ so that
$H$ corresponds to the purely real vectors in $\mathbb{O}^\ell$. Then
$x\cdot\mathfrak{h}\perp\mathfrak{h}$ if $x$ is a purely imaginary octonian and we obtain
a Riemannian submersion which is anti-invariant octonian.
\end{enumerate}
\item Let $m=2\ell$ and $n=\ell$. Identify $G$ with $\tilde{\mathbb{C}}^\ell$ so that
$H=\mathbb{R}^\ell$ corresponds to the purely real para-complex vectors.
More specifically, we take a basis $\{e_i,f_i\}$ for $\mathbb{R}^{2\ell}$ where
$H=\operatorname{Span}\{e_i\}$. Set
$$
\langle e_i,e_i\rangle=1,\ \langle f_i,f_i\rangle=-1,\ \langle e_i,f_j\rangle=0,\ \tilde Je_i=f_i,\ \tilde Jf_i=e_i\,.
$$
We obtain a Riemannian submersion which is Lagrangian para-Hermitian. By taking a different
inner product $\langle e_i,e_i\rangle=-\langle f_i,f_i\rangle=\epsilon_i$ for $\epsilon_i=\pm1$,
we can ensure that the base has arbitrary signature.
\item Let $m=4\ell$ and $n=\ell$. Identify $G=\tilde{\mathbb{Q}}^\ell$
so that $H=\mathbb{R}^\ell$ corresponds to
the purely real vectors in $G$. We
obtain a Riemannian submersion which is anti-invariant para-quaternion.
\end{enumerate}
\end{example}
The total space $G=\mathbb{R}^m$ is non-compact in Example~\ref{Ex3.1}.
We compactify by dividing by an integer lattice.

\begin{example}\label{Ex3.2}\rm Let $\mathbb{Z}^k$ be the integer lattice in $\mathbb{R}^k$
and let $\mathbb{T}^k:=\mathbb{R}^k/\mathbb{Z}^k$ be the
$k$-dimensional torus $S^1\times\dots\times S^1$ with the flat product metric. Let
$G=\mathbb{T}^m$ and let $H=\mathbb{T}^n$. We can
repeat the construction of Example~\ref{Ex3.1} to obtain examples which are compact.
\end{example}

\subsection{The Hopf fibration}
Example~\ref{Ex3.1} and Example~\ref{Ex3.2} are flat.
We can use the Hopf fibration to construct examples which are not flat.
We identify $\mathbb{R}^4$ with the quaternions $\mathbb{Q}$; this identifies
$S^3$ with the unit quaternions and gives $S^3$ a Lie group structure.
Let
\begin{equation}\label{Eq3.a}
e_1(x)=x\cdot i,\qquad e_2(x)=x\cdot j,\qquad e_3(x)=x\cdot k\,.
\end{equation}
This is then
a basis for the Lie algebra $\mathfrak{g}$ of left-invariant vector fields on $S^3$ and
$$[e_1,e_2]=-2e_3,\quad [e_2,e_3]=-2e_1,\quad [e_3,e_1]=-2e_2\,.$$
Every 1-dimensional Lie subalgebra of $S^3$ corresponds to a 1-dimensional compact
Abelian subgroup $S^1$ of $S^3$. Let $G=S^1\times S^3$ and let $e_0$ generate the Lie algebra of $S^1$ so that
$\mathfrak{g}=\operatorname{Span}\{e_0,e_1,e_2,e_3\}$. If $\epsilon\ne0$, define
\begin{equation}\label{Eq3.b}
\langle e_i,e_j\rangle=\left\{\begin{array}{rll}
1&\text{if }i=j=2\\
1&\text{if }i=j=3\\
\epsilon&\text{if }i=j=0\\
\epsilon&\text{if }i=j=1\\
0&\text{otherwise}\end{array}\right\}\,.
\end{equation}
These metrics are among the metrics first introduced by Hitchin~\cite{H74} in his
study of harmonic spinors and are Kaluza-Klein metrics.
Let
\begin{equation}\label{Eq3.c}\begin{array}{llll}
Je_0=e_1,&Je_1=-e_0,&Je_2=e_3,&Je_3=-e_2,\\
\check Je_0=e_2,&\check Je_2=-e_0,&\check Je_3=e_1,&\check Je_1=-e_3,\\
\tilde Je_0=e_2,&\tilde Je_2=e_0,&\tilde Je_1=-e_3,&\tilde Je_3=-e_1.
\end{array}\end{equation}
Then $J$ is an integrable Hermitian complex structure on $S^1\times S^3$ for any $\epsilon\ne0$. If $\epsilon=1$,
then $\{1,J,\check J,J\check J\}$ is a Hermitian quaternion structure on $S^1\times S^3$.
If $\epsilon=-1$, then $\tilde J$ is a para-Hermitian para-complex structure
and $\{1,J,\tilde J,J\tilde J\}$ is a Hermitian para-quaternion structure on $S^1\times S^3$.
If $\epsilon=1$, then $\langle\cdot,\cdot\rangle$ is bi-invariant. If $\epsilon\ne1$, then
$\langle\cdot,\cdot\rangle$ is right invariant under the 2-dimensional Lie subgroup $H$
with $\mathfrak{h}=\operatorname{Span}\{e_0,e_1\}$ but is not bi-invariant.

\begin{example}\label{Ex3.3}\rm
Let $G=S^1\times S^3$. Let $\mathfrak{h}$ be the Lie sub-algebra of a closed subgroup $H$ of $G$.
Adopt the notation of Equation~(\ref{Eq3.b}) and Equation~(\ref{Eq3.c}).
\begin{enumerate}
\item  Let $\mathfrak{h}=\operatorname{Span}\{e_0,e_1\}$.
\begin{enumerate}
\item Let $\epsilon=1$.
Then $G/H=S^2$ is the sphere of radius $2$ in $\mathbb{R}^3$ and has constant sectional curvature $\frac14$.
We use $\check J$ to obtain a Lagrangian Hermitian Riemannian submersion.
The
fibers of the submersion are minimal, not totally geodesic, and the horizontal
distribution is not integrable (see Park~\cite{P90}).
\item Let $\epsilon=-1$.
Then $G/H=S^2$. We use $\tilde J$ to
obtain a Lagrangian para-Hermitian Riemannian submersion.
\end{enumerate}
\item If $\mathfrak{h}=\operatorname{Span}\{e_0\}$, set $B=S^3$. If $\mathfrak{h}=\operatorname{Span}(e_1)$,
set $B=S^1\times S^3$.
\begin{enumerate}
\item Let $\epsilon\ne0$ be arbitrary.
We use $J$ to obtain
an anti-invariant Hermitian Riemannian submersion from $G$ to $B$.
\item Let $\epsilon=+1$. We use $J$ and $\check J$ to identify $\mathfrak{g}=\mathbb{Q}$ with the quaternions
to obtain an anti-invariant
quaternion Riemannian submersion from $G$ to $B$.
\item Let $\epsilon=-1$. We use $\tilde J$ to obtain
an anti-invariant Hermitian Riemannian submersion from $G$ to $B$.
\item Let $\epsilon=-1$. We use $J$ and $\tilde J$ to identify $\mathfrak{g}=\mathbb{Q}$ with the para-quaternions
to obtain an anti-invariant para-quaternion Riemannian
submersion from $G$ to $B$.
\end{enumerate}\end{enumerate}\end{example}

\subsection{Einstein geometry}

\begin{example}\label{Ex3.4}\rm
Let $g_{S^3}$ be the standard round metric on $S^3$ defined by $\langle e_i,e_j\rangle=\delta_{ij}$.
Let $G=(S^3)^\nu=S^3\times\dots\times S^3$. We take a product metric on $G$ where the metric on
each factor is $\pm g_{S^3}$;
thus this metric is bi-invariant. Let $H$ be a closed subgroup of $G$
and let $\pi:G\rightarrow G/H$ be
the associated Riemannian submersion.
\begin{enumerate}
\item Let $G=S^3\times S^3$ and $\mathfrak{g}=\operatorname{Span}\{e_1,e_2,e_3,f_1,f_2,f_3\}$.
\begin{enumerate}
\item
Let $g_G=g_{S^3}\oplus g_{S^3}$ be the
standard bi-invariant
Einstein metric on $S^3\times S^3$.
\begin{enumerate}
\item Let $\mathfrak{h}=\operatorname{Span}\{e_2,f_2\}$. Let $Je_1=f_1$, $Jf_1=-e_1$, $Je_2=e_3$, $Je_3=-e_2$, $Jf_2=f_3$, $Jf_3=-f_2$.
This almost complex structure is integrable and using $J$ yields
an anti-invariant Hermitian Riemannian submersion
from $S^3\times S^3$ to $S^2\times S^2$.
\item Let $\mathfrak{h}=\operatorname{Span}\{e_1,e_2,e_3\}$. Let $Je_i=f_i$ and $Jf_i=-e_i$
for $1\le i\le 3$. This almost
complex structure is not integrable. Using $J$ yields an anti-invariant almost Hermitian
Riemannian submersion from $S^3\times S^3$ to $S^3$
\end{enumerate}
\item Let $g_G=g_{S^3}\oplus -g_{S^3}$ be the standard
bi-invariant neutral signature metric on $S^3\times S^3$.
\begin{enumerate}
\item Let $\mathfrak{h}=\operatorname{Span}\{e_1,f_1\}$.
Let $\tilde Je_1=f_2$, $\tilde Jf_2=e_1$, $\tilde Jf_1=e_2$, $\tilde Je_2=f_1$, $\tilde Je_3=f_3$, and
$\tilde Jf_3=e_3$. Using $\tilde J$ yields an anti-invariant para-Hermitian
Riemannian submersion from $S^3\times S^3$ to $S^2\times S^2$.
 \item  Let $\mathfrak{h}=\operatorname{Span}\{e_1,e_2,e_3\}$.
 Let
$\tilde Je_i=f_i$ and $\tilde Jf_i=e_i$. Using $\tilde J$ yields
 an anti-invariant para-Hermitian Riemannian
submersion from $S^3\times S^3$ to $S^3$.
\end{enumerate}\end{enumerate}
\item Let $G=(S^3)^4$, let $g_G=g_{S^3}\oplus g_{S^3}\oplus g_{S^3}\oplus g_{S^3}$, and let $\dim\{H\}\le3$.
\begin{enumerate}
\item Identify $\mathfrak{g}$ with $\mathbb{Q}^3$ in such a way that
$\mathfrak{h}$ is real and the action of $\mathbb{Q}$ is Hermitian. Then $\pi$ is an anti-invariant quaternion.
\item Identify $\mathfrak{g}$ with $\tilde{\mathbb{Q}}^3$ in such a way
that $\mathfrak{h}$ is real and the action of $\tilde{\mathbb{Q}}$ is para-Hermitian. Then $\pi$
is an anti-invariant para-quaternion Riemannian submersion.
\end{enumerate}
\item Let $G=(S^3)^8$, let $g_G=g_{S^3}\oplus\dots g_{S^3}$, and let $\dim\{H\}\le 7$. Identify $\mathfrak{g}$ with
$\mathbb{O}^3$ in such a way that $\mathfrak{h}$ is real and the action of $\mathbb{O}$
is Hermitian. Then $\pi$ is an anti-invariant octonian Riemannian submersion.
\end{enumerate}\end{example}

\subsection{Negative curvature}
Our previous examples have, for the most part, involved the Lie group $S^3$ and
the Hopf fibration $S^3\rightarrow S^2$. We now turn to the negative curvature dual.
Let $\mathbb{H}^2(0,2)$ be the hyperbolic plane with a Riemannian metric
of constant sectional curvature $-\frac14$ and let $\mathbb{H}^2(1,1)$ be the Lorentzian analogue.
We recall some facts about the
Lie group $\operatorname{SL}(2,\mathbb{R})$ and refer to Section~6.8 of
Gilkey, Park, and V\'azquez-Lorenzo \cite{GPV15} -- there are, of course, many excellent references.
$\operatorname{SL}(2,\mathbb{R})$ is a $3$-dimensional Lie group
and the Lie algebra $\mathfrak{sl}(2,\mathbb{R})$ is the vector space of trace free
$2\times 2$ real matrices. The canonical basis for
$\mathfrak{sl}(2,\mathbb{R})$ is
$$f_1:=\left(\begin{array}{rr}0&1\\-1&0\end{array}\right),\quad
 f_2:=\left(\begin{array}{rr}0&1\\1&0\end{array}\right),\quad
 f_3:=\left(\begin{array}{rr}1&0\\0&-1\end{array}\right).
$$
The bracket relations then take the form
$$[f_1,f_2]=2f_3,\quad[f_2,f_3]=-2f_1,\quad[f_3,f_1]=2f_2\,.$$
The Lie algebra $\mathfrak{s}^3$ of $S^3$ is the Lie algebra of the special unitary group $SU(2)$ in positive
definite signature and
the Lie algebra of $\mathfrak{sl}(2,\mathbb{R})$ of $\operatorname{SL}(2,\mathbb{R})$
is the Lie algebra of the special unitary group
$SU(1,1)$ in indefinite signature; the two
are related by complexification.
Let $\operatorname{ad}(\xi):\eta\rightarrow[\xi,\eta]$
be the adjoint action and let
$K(\xi,\eta):=\operatorname{Tr}\{\operatorname{ad}(\xi)\operatorname{ad}(\eta)\}$
be the Killing form. One then has
$$K(f_i,f_j)=\left\{\begin{array}{l}
-8\text{ if }i=j=1\\+8\text{ if }i=j=2\\
+8\text{ if }i=j=3\\
\phantom{+}0\text{ otherwise}\end{array}\right\}\,.
$$
There is no bi-invariant Riemannian metric on $\operatorname{SL}(2,\mathbb{R})$.
However, $\frac18K$ is a bi-invariant
Lorentzian metric on $\operatorname{SL}(2,\mathbb{R})$. Let
\begin{eqnarray*}
&&\sigma_1(x):=\left(\begin{array}{rr}\cos(x)&\sin(x)\\-\sin(x)&\cos(x)\end{array}\right),\quad
\sigma_2(x):=\left(\begin{array}{rr}\cosh(x)&\sinh(x)\\\sinh(x)&\cosh(x)\end{array}\right),\\
&&\sigma_3(x):=\left(\begin{array}{rr}e^x&0\\0&e^{-x}\end{array}\right)\,.
\end{eqnarray*}
These define closed Abelian Lie sub-groups $H_i$ of $\operatorname{SL}(2,\mathbb{R})$
whose associated Lie-algebras are spanned by $f_i$. The natural coset spaces
$\operatorname{SL}(2,\mathbb{R})/H_i$ have constant negative sectional curvature $-\frac14$
and may be identified with $\mathbb{H}^2(0,2)$ if $i=1$ and $\mathbb{H}^2(1,1)$ if $i=2,3$.

Let $G=\mathbb{R}\times\operatorname{SL}(2,\mathbb{R})$. Let $f_0$ correspond to the Abelian factor.
Define a bi-invariant neutral signature metric on $G$ by setting:
$$
\langle f_i,f_j\rangle=\left\{\begin{array}{l}-1\text{ if }i=j=0\\
-1\text{ if }i=j=1\\
+1\text{ if }i=j=2\\
+1\text{ if }i=j=3\\
\phantom{+}0\text{ otherwise}
\end{array}\right\}\,,
$$
In analogy with Equation~(\ref{Eq3.b}), we set:
$$\begin{array}{llll}
Jf_0=f_1,&Jf_1=-f_0,&Jf_2=f_3,&Jf_3=-f_2,\\
\tilde Jf_0=f_2,&\tilde Jf_2=f_0,&\tilde Jf_1=-f_3,&\tilde Jf_3=-f_1.
\end{array}$$
Then $J$ is a Hermitian complex structure on $G$ and $\tilde J$ is a Hermitian para-complex
structure on $G$; $J$ and $\tilde J$ generate a Hermitian para-complex structure on $G$.

\begin{example}\label{E3.6}\rm
Let $G=\mathbb{R}\times\operatorname{SL}(2,\mathbb{R})$,
let $H$ be a closed subgroup of $G$, and let $\pi$ be the natural projection from $G$ to $G/H$.
\begin{enumerate}
\item If $\mathfrak{h}=\operatorname{Span}\{f_0\}$, then $\pi$ is an anti-invariant Hermitian,
an anti-invariant para-Hermitian, and
an anti-invariant para-quaternion Riemannian submersion from $G$ to $\operatorname{SL}(2,\mathbb{R})$.
\item  If $\mathfrak{h}=\operatorname{Span}\{f_1\}$, then $\pi$ is
an anti-invariant Hermitian, an anti-invariant para-Hermitian, and
an anti-invariant para-quaternion Riemannian submersion from $G$ to $\mathbb{R}\times\mathbb{H}^2(0,2)$.
\item If $\mathfrak{h}=\operatorname{Span}\{f_2\}$, then $\pi$ is
an anti-invariant Hermitian, an anti-invariant para-Hermitian, and
an anti-invariant para-quaternion Riemannian submersion from $G$ to $\mathbb{R}\times\mathbb{R}\mathbb{H}^2(1,1)$.
\item If $\mathfrak{h}=\operatorname{Span}\{f_0,f_2\}$,
then $\pi:G\rightarrow G/H$ is an anti-invariant Hermitian Riemannian submersion from $G$ to $\mathbb{H}^2(1,1)$.
\end{enumerate}
\end{example}

\section{Examples where the total space is not a Lie group}\label{S4}
 In this section, we present examples where the total space is
not a Lie group. We identify $\mathbb{R}^{4k}$ with
$\mathbb{C}^{2k}$ to define an action of $S^1$ on $S^{4k-1}$; the
quotient $S^{4k-1}/S^1$ is complex projective space
$\mathbb{CP}^{2k-1}$ with a Fubini-Study metric of constant positive
holomorphic sectional curvature. We identify $\mathbb{R}^{4\ell}$
with $\mathbb{Q}^\ell$ to define an action of $S^3$ on
$S^{4\ell-1}$; the quotient $S^{4\ell-1}/S^3$ is quaternionic
projective space $\mathbb{QP}^{\ell-1}$. Instead of taking the
Euclidean inner product on $\mathbb{R}^{2k}$, we could take an
indefinite signature metric. Let
\begin{equation}\begin{array}{l}\label{E4.a}
\langle x,y\rangle:=-x^1y^1-x^2y^2+x^3y^3+\dots+x^{2\ell}y^{2\ell},\\[0.05in]
\ll x,y\gg:=-x^1y^1-\dots-x^4y^4+x^5y^5+\dots x^{4\ell}y^{4\ell}\\[0.05in]
\tilde S^{2k-1}:=\langle x,y\rangle=-1,\quad\check S^{4k-1}:=\{x\in\mathbb{R}^{4\ell}:\ll x,x\gg=-1\}\,.
\end{array}\end{equation}
The pseudo-spheres $\tilde S^{2k-1}$ and $\check S^{4k-1}$
inherit indefinite signature metrics of constant sectional curvature.
The quotient $\tilde S^{2k-1}/S^1$ is the negative curvature dual of $\mathbb{CP}^{k-1}$ and the quotient
$\check S^{4k-1}/S^3$ is the negative
curvature of $\mathbb{QP}^{k-1}$.

Let $m=2k$ and $N=S^{2k-1}$ or $N=\tilde S^{2k-1}$, or let $m=4k$ and $N=\check S^{4k-1}$.
There is an orthogonal direct sum decomposition
$T(\mathbb{R}^m)|_N=\nu\oplus T(N)$
where $\nu$ is the normal bundle. Let $M=S^1\times N$ with the product metric.
Since $\nu$ is a trivial line bundle, we have a natural isometry
$\Xi:TM\approx M\times\mathbb{R}^n$. Let $\partial_\theta$ be the natural unit tangent vector field on $S^1$.
Let $x=(\theta,\Theta)\in M$. Then $\Xi(\theta,\Theta)\nu=\Theta$ and $\Xi(TN)=\Theta^\perp$.

\begin{example}\label{E4.1}\rm
Let $\ell\ge3$. {Adopt the notation of Equation~(\ref{E4.a})}. Let
$J$ be complex multiplication by $i$ on $TM=M\times\mathbb{C}^\ell$.
Let $H=S^1$. Since $\dim\{\mathcal{V}\}=1$,
$J\mathcal{V}\perp\mathcal{V}$.
\begin{enumerate}
\item Let $M=S^1\times S^{2\ell-1}$. Let $H$ act on $S^1$ by complex multiplication and trivially on $S^{2\ell-1}$.
Then $\pi:M\rightarrow S^{2\ell-1}$
is an anti-invariant Hermitian Riemannian submersion.
\item Let $M=S^1\times S^{2\ell-1}$. Let $H$ act on trivially on $S^1$ and by complex multiplication on $S^{2\ell-1}$.
Then $\pi:M\rightarrow S^1\times\mathbb{CP}^{\ell-1}$
is an anti-invariant Hermitian Riemannian submersion.
\item Let $M=S^1\times\tilde S^{2\ell-1}$. Let $H$ act on $S^1$ by complex multiplication and trivially on $\tilde S^{2\ell-1}$.
Then $\pi:M\rightarrow\tilde S^{2\ell-1}$
is an anti-invariant Hermitian Riemannian submersion.
\item Let $M=S^1\times\tilde S^{2\ell-1}$. Let $H$ act on trivially on $S^1$ and by complex multiplication on $\tilde S^{2\ell-1}$.
Then
$\pi:M\rightarrow S^1\times\widetilde{\mathbb{CP}}^{\ell-1}$
is an anti-invariant Hermitian Riemannian submersion.
\end{enumerate}
We have taken $\ell\ge3$ since the case $\ell=2$ recovers the Hopf
fibration $S^3\rightarrow S^2$ or $S^3\rightarrow\mathbb{H}^2$.
\end{example}

\begin{example}\label{Ex4.3}\rm {Adopt the notation of Equation~(\ref{E4.a})}.
Let $\ell\ge2$. Let $H=S^1\times S^1$. Let $J$ be quaternion multiplication on $TM=M\times\mathbb{Q}^\ell$.
\begin{enumerate}
\item Use the product action to let $H$ act on the first and on the second factor
of $M=S^1\times S^{4\ell-1}$. Let $\pi$ be the associated Riemannian submersion
from $M$ to $\mathbb{CP}^{2\ell-1}$.
Then $\mathcal{V}(\theta,\Theta)=\operatorname{Span}\{\Theta,i\cdot\Theta\}$.
Since $j\cdot\mathcal{V}\perp\mathcal{V}$, $\pi$ is an anti-invariant Hermitian Riemannian submersion.
\item Use the product action to let $H$ act on the first and on the second factor
of $M=S^1\times\check S^{4\ell-1}$. Let $\pi$ be the associated Riemannian submersion
from $M$ to $\widetilde{\mathbb{CP}}{}^{2\ell-1}$.
Then $\mathcal{V}(\theta,\Theta)=\operatorname{Span}\{\Theta,i\cdot\Theta\}$.
Since $j\cdot\mathcal{V}\perp\mathcal{V}$, $\pi$ is an anti-invariant Hermitian Riemannian submersion.
\end{enumerate}
We have taken $\ell\ge2$ since the case $\ell=1$ recovers the Hopf
fibration $S^3\rightarrow S^2$ or $S^3\rightarrow\mathbb{H}^2$.
\end{example}

\section*{Acknowledgments}
\noindent Research of the authors was partially supported by the
Basic Science Research Program through the National Research
Foundation of Korea(NRF) funded by the Ministry of Education
(2014053413). {It is a pleasant task to acknowledge helpful
correspondence concerning these matters with Professor B.
{\c{S}ahin}}.


\begin{thebibliography}{lll}

\bibitem{AF13} S. Ali and T. Fatima, ``Anti-invariant Riemannian submersions from
nearly K\"ahler manifolds", {\it Filomat \bf 27} (2013), 1219--1235. DOI 10.2298/FIL.1307219A.

\bibitem{AF13a} S. Ali and T. Fatima, ``Generic Riemannian submersions",
{\it Tamkang J. of Math. \bf 44} (2013), 395--405.

\bibitem{AM01} D. V. Alekseevsky and S. Marchiafava,
``Almost complex submanifolds of quaternionic manifolds", {\it
 Proceedings of the colloquium on differential geometry, Debrecen (Hungary)} , 25-30 July 2000, Inst. Math. Inform. Debrecen, 2001, 23--38.

\bibitem{A13} M. Atceken, ``Anti-invariant Riemannian submersions from locally Riemannian
product manifold to a Riemannian manifold", {\it Gulf J. Mathematics \bf 1} (2013), 25--35.

\bibitem{BKEM15} A. Beri, I. K\"upeli Erken, and C. Murathan,
``Anti-invariant Riemannian submersions from Kenmotsu manifolds onto Riemannian manifolds",
arXiv:1504.04180.

\bibitem{GPV15} P. Gilkey, J. H. Park and R. V\'azquez-Lorenzo, ``Aspects of Differential Geometry II",
Morgan and Claypool (2015), doi:10.2200/S00645ED1V01Y201505MAS016.

\bibitem{G13} Y. G\"und\"uzalp, ``Anti-invariant semi-Riemannian submersions from
almost para-Hermitian manifolds", J. Function Spaces and Applications (2013), 7020623.

\bibitem{G14} Y. G\"und\"uzalp, ``Anti-invariant semi-Riemannian submersions from
almost product Riemannian manifolds",  Mathematical Sciences And
Applications E-Notes Volume 1 No. 1 (2013), 58-66.

\bibitem{H74} N. Hitchin, ``Harmonic Spinors", {\it Adv. in Math. \bf 14} (1974), 1--65.

\bibitem{IZ11} S. Ivanov and S. Zamkovoy, ``Parahermition and paraquaternionic manifolds",
{\it J. Diff. Geo. and Appl. \bf 23} (2005), 206--234.


\bibitem{LPSS15} J. Lee, J. H. Park, B. {\c{S}ahin}, and D.-Y. Song,
``Einstein conditions for the base space of anti-invariant Riemannian submersions and Clairaut
submersions", {\it Taiwanese J. Math \bf 19}, 1145--1160.

\bibitem{MK13} C. Murthan and I. K\"upeli Erkin, ``Anti-invariant Riemannian
submersions from cosymplectic manifolds", {Filomat 29 (7), (2015),
1429--1444.}

\bibitem{NN57} A. Newlander and L. Nirenberg,
``Complex analytic coordinates in almost complex manifolds",
{\it Annals of Mathematics \bf65} (1957), 391--404. doi:10.2307/1970051

\bibitem{O66} B. O'Neill, ``The fundamental equations of a submersion",
{\it Mich. Math. J. \bf13} (1966), 458--469.

\bibitem{P90} J. H. Park, ``The Laplace-Beltrami operator and
Romanian submersion with minimal and not totally geodesic fibers",
{\it Bull. Korean Math. Soc. \bf 27} (1990), 39--47.

\bibitem{P15} K. Park, ``H-anti-invariant submersions from almost quaternionic Hermitian manifolds",
arXiv 1507.04473v1 [math.DG].

\bibitem{S10} B. {\c{S}ahin}, ``Anti-invariant Riemannian submersions from almost Hermitian
manifolds", {\it Central European J. Math \bf 8} (2010), 437--447.
DOI 10.2478/s11533-010-0023-6,

\bibitem{S13} B. {\c{S}ahin}, ``Riemannian submersions from almost Hermitian Manifolds",
{\it Taiwanese Journal of Mathematics \bf17} (2013), 629--659.

\bibitem{W15} Wikipedia, ``https://en.wikipedia.org/wiki/octonian".

\end{thebibliography}
\end{document}